\documentclass[12pt]{amsart}
\catcode`\@=11 \@addtoreset{equation}{section}

\catcode`\@=12

\newtheorem{Theorem}{Theorem}[section]
\newtheorem{Proposition}{Proposition}[section]
\newtheorem{Lemma}{Lemma}[section]
\newtheorem{Corollary}{Corollary}[section]

\newcommand{\bTheorem}[1]{
\begin{Theorem} \label{T#1} }
\newcommand{\eT}{\end{Theorem}}

\newcommand{\bProposition}[1]{
\begin{Proposition} \label{P#1}}

\newcommand{\eP}{\end{Proposition}}

\newcommand{\bLemma}[1]{
\begin{Lemma} \label{L#1} }
\newcommand{\eL}{\end{Lemma}}

\newcommand{\bCorollary}[1]{
\begin{Corollary} \label{C#1} }
\newcommand{\eC}{\end{Corollary}}

\newcommand{\bFormula}[1]{
\begin{equation} \label{#1}}
\newcommand{\eF}{\end{equation}}

\newcommand{\Ov}[1]{\overline{#1}}
\newcommand{\VC}[2]{ \left[ \begin{array}{c} #1 \\ #2 \end{array} \right] }
\newcommand{\DC}{C^\infty_c}

\newcommand{\vr}{\varrho}
\newcommand{\vre}{\vr_\ep}

\newcommand{\vue}{\vu_\ep}

\newcommand{\vu}{\vc{u}}
\newcommand{\vc}[1]{{\mathbf {#1}}}
\newcommand{\Div}{{\rm div}_x}
\newcommand{\Grad}{\nabla_x}

\newcommand{\tn}[1]{\mbox {\F #1}}
\newcommand{\dx}{{\rm d} {x}}
\newcommand{\dt}{{\rm d} t }

\newcommand{\Rm}{\mbox{\FF R}}

\newcommand{\intO}[1]{\int_{\Omega} #1 \ \dx}

\newcommand{\ep}{\varepsilon}

\font\F=msbm10 scaled 1100
\font\FF=msbm10 scaled 800

\font\FF=msbm10 scaled 800
\def\longformule#1#2{
\displaylines{ \qquad{#1} \hfill\cr \hfill {#2} \qquad\cr } }
\newcommand{\R}{\mathbb{R}}

\date{}
\begin{document}

\title{A singular limit for compressible rotating fluids}

\author{Eduard Feireisl}
\address[Eduard Feireisl]%
{Institute of Mathematics of the Academy of Sciences of the Czech Republic \\
\v Zitn\' a 25, 115 67 Praha 1, Czech Republic}
\email{feireisl@math.cas.cz}
\thanks{The work of E.F. was supported by Grant 201/08/0315 of GA \v CR as a part of the general research
programme of the Academy of Sciences of the Czech Republic,
Institutional Research Plan AV0Z10190503.}

\author{Isabelle Gallagher}
\thanks{The work of I.G is partially supported by the French Ministry of Research grant
ANR-08-BLAN-0301-01.}
\address[Isabelle Gallagher]%
{Institut de Math{\'e}matiques UMR 7586 \\
      Universit{\'e} Paris Diderot \\
175, rue du Chevaleret\\
75013 Paris\\ France}
\email{Isabelle.Gallagher@math.jussieu.fr}
\author{Anton{\' \i}n Novotn\' y}
\address[Anton{\' \i}n Novotn\' y]%
{IMATH Universit\'e du Sud Toulon-Var \\
B.P. 132\\
83957 La Garde Cedex\\
France}
\email{novotny@univ-tln.fr}
\maketitle
%
%\centerline{Institute of Mathematics of the Academy of Sciences of the Czech Republic}
%\centerline{\v Zitn\' a 25, 115 67 Praha 1, Czech Republic}

%
%\centerline{Institut de Math\'ematiques de Jussieu and Universit\'e Paris Diderot}
%\centerline{175 rue du Chevaleret, 75013 Paris, France}

\begin{abstract}
We consider a singular limit problem for the Navier-Stokes system of a rotating compressible fluid, where the Rossby and Mach numbers tend simultaneously to zero. The limit problem is identified as the 2-D Navier-Stokes system in the ``horizontal'' variables containing an extra term that accounts for compressibility in the original system.
\end{abstract}
\section{Introduction}
\label{i}
Consider a scaled Navier-Stokes system in the form
\bFormula{i1}
\partial_t \vr + \Div (\vr \vu) = 0,
\eF
\bFormula{i2}
\partial_t (\vr \vu) + \Div (\vr \vu \otimes \vu) + \frac{1}{\ep} ( {\vc{g}} \times
\vr \vu ) + \frac{1}{\ep^2} \Grad p(\vr) =  \Div \tn{S}(\Grad \vu ),
\eF
with the viscous stress tensor
\bFormula{i3}
\tn{S}(\Grad \vu)  = \mu \Big( \Grad \vu + \Grad^t \vu - \frac{2}{3} \Div \vu \tn{I} \Big),\
\mu > 0,
\eF
and
\[
\vc{g} = [0,0,1].
\]
Here $ \vr=\vr(t,x)\ge 0$ denotes the  density and $\vc u(t,x)=[u_1,u_2,u_3](t,x)$ denotes
the velocity of the fluid.
Problem (\ref{i1} - \ref{i2}) arises in meteorological applications, modeling rotating compressible fluids with the rotation axis determined by $\vc{g}$ and the Rossby and Mach number proportional to a small parameter $\ep$.

We consider a very simple geometry of the underlying physical space, namely
an infinite slab $\Omega$ bounded above and below by two parallel planes,
\bFormula{i4}
\Omega = \R^2 \times (0,1).
\eF
The velocity $\vu$ satisfies the complete slip boundary conditions,
\bFormula{i5}
\vu \cdot \vc{n} = u_3 |_{\partial \Omega} = 0,\
[\tn{S} \vc{n}] \times \vc{n}|_{\partial \Omega} = [S_{2,3}, - S_{1,3}, 0 ]|_{\partial \Omega} = 0.
\eF
For the initial data
\bFormula{i6}
\vr(0, \cdot) = \vr_{0, \ep}, \ \vc{u}(0,\cdot) = \vu_{0,\ep},
\eF
our goal is to study the asymptotic behavior of the corresponding solutions~$\vre$, $\vue$ for $\ep \to 0$. We focus on the interplay between the Coriolis force,
here proportional to a singular parameter $1/\ep$, and the acoustic waves created in the low Mach number regime. In particular, we neglect:

\begin{itemize}

\item stratification due to the presence of gravitation, here assumed in equilibrium with the centrifugal force; accordingly
the action of the centrifugal force is also neglected;

\item the effect of a boundary layer (Ekman layer), here eliminated by the choice of the
complete slip boundary conditions.

\end{itemize}

We consider   \emph{ill-prepared} initial data, specifically,
\bFormula{i7}
\left\{
\begin{array}{c}
\vr_{0,\ep} = \Ov{\vr} + \ep r_{0,\ep}, \ \mbox{with}\ \{ r_{0,\ep} \}_{\ep > 0}
\ \mbox{bounded in} \ L^2 \cap L^\infty (\Omega),\\ \\
\mbox{for some positive constant}\ \Ov{\vr}, \\ \\
\{ \vu_{0, \ep} \}_{\ep > 0} \ \mbox{bounded in} \ L^2 \cap L^\infty(\Omega;\R^3).
\end{array} \right\}
\eF

Because of the prominent role of the ``vertical'' direction $\vc{g}$ in the problem, we introduce the ``horizontal'' component $\vc{v}_h = [v_1,v_2,0]$ of a vector field $\vc{v}$, together with the corresponding differential operators~$\nabla_h$, ${\rm div}_h$, and, notably, ${\rm curl}_h$, which is represented by the \emph{scalar} field
\[
{\rm curl}_h [\vc{v}] = \partial_{x_1} v_2 - \partial_{x_2} v_1.
\]
Let $\vre$, $\vue$ be a solution of problem (\ref{i1}-\ref{i6}). Introducing a new quantity
\[
r_\ep = \frac{ \vre - \Ov{\vr} }{\ep}
\]
which satisfies
\[
\partial_t r_\ep +\frac 1 \ep \Div ( \Ov{\vr} \vu)  +  \Div (  r_\ep \vu)  = 0,
\]
we easily check that if
\[
r_\ep \to r,\ \vue \to \vc{U} \ \mbox{in some sense,}
\]
then, at least formally, the limits satisfy a diagnostic equation
\bFormula{i8}
\vc{g} \times \vc{U} + \frac{p'(\Ov{\vr})}{\Ov{\vr}} \Grad r = 0,
\eF
which in turn implies that
\bFormula{i9}
r = r(x_1,x_2), \ \vc{U} = [\vc{U}_h,0], \ \vc{U}_h = \vc{U}_h (x_1,x_2).
\eF
Moreover, as we will see below, $\Div \vc{U} = {\rm div}_h \vc{U}_h = 0$, and denoting~$ \nabla_h^\perp r$ the vector~$(\partial_{x_2}r, -\partial_{x_1}r)$,
\bFormula{i10}
\partial_t \Big( \Delta_h r - \frac{1}{p'(\Ov{\vr})} r \Big) +
\nabla^\perp_h r \cdot \nabla_h ( \Delta_h r ) = \frac{\mu}{\Ov{\vr}} \Delta^2_h r.
\eF
Thus $r$ may be interpreted as a \emph{stream function} associated to the vector field
$\vc{U}_h$, therefore (\ref{i10}) can be viewed as a 2D Navier-Stokes system describing the motion of an incompressible fluid in the horizontal plane~$\R^2$, supplemented with an extra term $(1/p'(\Ov{\vr})) \partial_t r$.

The main goal of the present paper is to provide a rigorous justification of the target system (\ref{i10}) in the framework of weak solutions to the primitive equations
(\ref{i1}), (\ref{i2}). In Section \ref{w}, we introduce the weak solutions to both systems,  recall their basic properties, and state our main result. In Section \ref{a}, we derive the necessary uniform bounds on the family
of solutions $\{ \vre, \vue \}_{\ep > 0}$, and pass formally to the limit when $\ep \to 0$.
In Section \ref{aw}, the associated wave equation describing propagation of the acoustic waves in the low Mach number regime is introduced. Using the celebrated RAGE theorem, we show that the acoustic energy tends to zero, at least locally in space. The proof of the main result is completed in Section \ref{c}.

 \section{Preliminaries}
\label{w}

To begin, we point out that system (\ref{i1} - \ref{i3}), endowed with the boundary conditions (\ref{i5}) can be recast as a purely periodic problem with respect to the vertical coordinate $x_3$ provided $\vr$, $u_1$, $u_2$ were extended as even functions
in the $x_3-$variable defined on
\[
\Omega = \R^2 \times {\mathcal T}^1, \ {\mathcal T}^1 \equiv [-1,1]|_{\{ -1,1 \}},
\]
while $u_3$ is extended to be odd in $x_3$ on the same set. A similar convention is adopted for the initial data.

\subsection{Weak solutions}

We shall say that functions $\vr$, $\vu$ represent a \emph{weak solution} to problem
(\ref{i1} - \ref{i6}) in $(0,T) \times \Omega$ if:
\begin{itemize}
\item
$\vr \geq 0$, $\vr \in L^\infty(0,T;L^\gamma(\Omega))$ for a certain $\gamma > 3/2$,

$\vu \in L^2(0,T; W^{1,2}(\Omega;\R^3))$;
\item
equation of continuity (\ref{i1}) is satisfied in the sense of renormalized solutions, namely
\bFormula{i11}
\int_0^T \intO{ \Big( (\vr + b(\vr)) \partial_t \varphi + (\vr + b(\vr)) \vu \cdot \Grad \varphi +
(b(\vr) - b'(\vr) \vr) \Div \vu \varphi \Big) } \ \dt
\eF
\[
= - \intO{ \Big( \vr_{0,\ep} + b(\vr_{0, \ep} ) \Big) \varphi (0, \cdot) }
\]
for any $b \in C^\infty[0, \infty)$, $b' \in \DC[0, \infty)$, and any test function
$\varphi \in \DC([0,T) \times \Omega)$;

\item
$p = p(\vr) \in L^1((0,T) \times \Omega)$,
momentum equation (\ref{i2}) is replaced by a family of integral identities
\bFormula{i12}
\int_0^T \intO{ \Big( \vr \vu \cdot \partial_t \varphi + \vr (\vu \otimes \vu) :
\Grad \varphi + \frac{1}{\ep}(\vc{g} \times \vr \vu) \cdot \varphi + \frac{1}{\ep^2} p(\vr) \Div \varphi \Big)} \ \dt
\eF
\[
= \int_0^T \intO{ \tn{S}(\Grad \vu) : \Grad \varphi }\ \dt - \intO{ \vr_{0, \ep} \vu_{0, \ep}
\cdot \varphi(0, \cdot) }
\]
for any $\varphi \in \DC([0,T) \times \Omega; \R^3)$;
\item
the \emph{energy inequality}
\bFormula{i13}
\intO{ \left( \frac{1}{2} \vr |\vu|^2 + \frac{1}{\ep^2} E(\vr, \Ov \vr) \right) (\tau, \cdot) } +
\int_0^\tau \intO{ \tn{S}( \Grad \vu ) : \Grad \vu }  \ \dt
\eF
\[
\leq
\intO{ \left( \frac{1}{2} \vr_{0,\ep} |\vu_{0,\ep}|^2 + \frac{1}{\ep^2} E(\vr_{0,\ep}, \Ov \vr) \right) }
\]
holds for a.a. $\tau \in (0,T)$, where
\[
E(\vr, \Ov{\vr}) = H(\vr) - H'(\Ov{\vr})(\vr - \Ov{\vr}) - H(\Ov{\vr}),
\]
with
\[
H(\vr) = \vr \int_1^\vr \frac{p(z)}{z^2} \ {\rm d}z.
\]
\end{itemize}

Note that, by virtue of hypothesis (\ref{i7}), the quantity on the right-hand side of
(\ref{i13}) is bounded uniformly for $\ep \to 0$.

\emph{Existence} of global-in-time weak solutions to problem (\ref{i1} - \ref{i6}) can be established by the method developed by P.-L. Lions \cite{LI4}, with the necessary modifications introduced in \cite{FNP} in order to accommodate a larger class of physically relevant pressure-density state equations, specifically,
\bFormula{i14}
p \in C^1[0, \infty), \ p(0) = 0, \ p'(\vr) > 0 \ \mbox{for}\ \ \vr > 0,\
\lim_{\vr \to \infty} \frac{p'(\vr)}{\vr^{\gamma - 1}} = p_\infty > 0,
\eF
for a certain $\gamma > 3/2$.

\subsection{Main result}

The main result of the present paper can be stated as follows.

\bTheorem{w1}
Assume that the pressure $p$ satisfies (\ref{i14}).

Let~$\{ \vre, \vue \}_{\ep > 0}$ be a family of weak solutions to problem
(\ref{i1} - \ref{i6}) in $(0,T) \times \Omega$, where $\Omega$ is specified through (\ref{i4}), with the initial data satisfying (\ref{i7}), where
\[
r_{0,\ep} \to r_0 \ \mbox{weakly in}\ L^2(\Omega),\
\vu_{0,\ep} \to \vc{U}_0 \ \mbox{weakly in}\ L^2(\Omega; \R^3).
\]

Then after taking a subsequence,  the following results hold
\[
r_\ep \equiv \frac{ \vre - \Ov{\vr} }{\ep}  \to r \ \mbox{weakly-(*) in}\ L^\infty(0,T; L^2(\Omega) + L^\gamma(\Omega)),
\]
\[
\vue \to \vc{U} \ \mbox{weakly in} \ L^2(0,T; W^{1,2}(\Omega; \R^3)),
\]
\[
 \mbox{and}†\ \vue \to \vc{U} \ \mbox{strongly in}  \ L^2_{\rm loc}((0,T) \times \Omega; \R^3),
\]
where $r$ and $\vc{U}$ satisfy (\ref{i8}), $\Div \vc{U} = 0$, and, moreover, the stream function~$r$ solves equation
(\ref{i10}) in the sense of distributions, supplemented with the initial datum
\bFormula{idat}
r(0, \cdot) = \tilde r ,
\eF
where $\tilde r \in W^{1,2}(\R^2)$ is the unique solution of
\[
- \Delta_h \tilde r +\frac1{p'( \Ov{\vr})} \tilde r =  \Ov{\vr} \int_0^1 {\rm curl}_h \vc{U}_{0,h} \ {\rm d}x_3 +
\int_0^1 r_0 {\rm d}x_3.
\]

If, in addition, ${\rm curl_h} \vc{U}_{0,h} \in L^2(\Omega)$,
then the solution $r$ of (\ref{i10}) is uniquely determined by (\ref{idat}) and the convergence holds for the whole sequence of solutions.
\eT

The remaining part of the paper is devoted to the proof of Theorem~\ref{Tw1}.
The crucial point of the proof is, of course, the strong (a.a. pointwise) convergence of the velocity field that enables us to carry out the limit in the convective term. Here, the desired pointwise convergence will follow from the celebrated RAGE theorem, together with the fact that the wave propagator in the associated \emph{acoustic equation} commutes with the Fourier transform in both the horizontal variables~$(x_1,x_2)$ and the vertical variable $x_3$.

\subsection{Related results}
This work is a contribution to a general research direction consisting in studying singular limits in PDEs arising in fluid mechanics. Without giving an extensive bibliography, one should refer for the first works in this line to Klainerman and Majda \cite{klainermanmajda} and Ukai \cite{ukai} for the incompressible limit (actually~\cite{ukai} is probably the first work in which dispersive estimates were established in order to prove strong convergence in the whole space), followed  by Desjardins et al. \cite{desjardinsgrenier} and~\cite{dglm}. In the context of rotating fluids one should mention the important work of Babin, Mahalov and Nicolaenko \cite{bmn}, as well as the book~\cite{cdgg} and references therein; one also refers to~\cite{GSR2} for a survey. Few studies combine both rotation and compressible effects. We refer to Bresch, Desjardins and G\' erard-Varet \cite{BDG} for an analysis in a cylinder, where the well prepared case is studied precisely; the ill prepared case is also addressed but only a conditional result is proved.

\section{Uniform bounds}
\label{a}

We start reviewing rather standard uniform bounds that follow directly from the energy inequality (\ref{i13}). To this end, it is convenient to introduce a decomposition
\[
h = [h]_{\rm ess} + [h]_{\rm res}, \ \mbox{where}\ [h]_{\rm ess} = \psi (\vre) h,
\]
\[
\psi \in \DC(0,\infty), \ 0 \leq \psi \leq 1, \ \psi \equiv 1 \ \mbox{in a neighborhood of}
\ \Ov{\vr}
\]
for any function $h$ defined on $(0,T) \times \Omega$.
It is understood that the \emph{essential} part $[h]_{\rm ess}$ is the crucial quantity that determines the asymptotic behavior of the system while the \emph{residual} component
$[h]_{\rm res}$ ``disappears'' in the limit $\ep \to 0$.

As already pointed out, our choice of the initial data (\ref{i7}) guarantees that the right-hand side of energy inequality (\ref{i13}) remains bounded for~$\ep \to 0$. After a straightforward manipulation, we deduce the following estimates:
\bFormula{a1}
\{ \sqrt{\vre} \vue \}_{\ep > 0} \ \mbox{bounded in} \ L^\infty(0,T; L^2(\Omega; \R^3)),
\eF
\bFormula{a2}
\{ [ r_\ep ]_{\rm ess} \}_{\ep > 0} \ \mbox{bounded in} \ L^\infty(0,T; L^2(\Omega)),
\eF
\bFormula{a3}
{\rm ess} \sup_{t \in (0,T)} \| [\vre]_{\rm res} \|^\gamma_{L^\gamma(\Omega)} \leq \ep^2 c,
\eF
\bFormula{a4}
{\rm ess} \sup_{t \in (0,T)} \| [1]_{\rm res} \|_{L^1(\Omega)} \leq \ep^2,
\eF
and
\bFormula{a5}
\left\{ \Grad \vue + \Grad^t \vue - \frac{2}{3} \Div \vue \tn{I} \right\}_{\ep > 0}
\: \mbox{bounded in} \  L^2((0,T) \times \Omega; \R^{3 \times 3} ).
\eF

In addition, it is easy to observe that (\ref{a2}), (\ref{a3}) yield
\bFormula{a6}
\vre \to \Ov{\vr} \ \mbox{in} \ L^\infty(0,T; L^\gamma + L^2(\Omega)),
\eF
which, together with (\ref{a1}), (\ref{a6}) and the standard Korn  inequality, gives rise to
\bFormula{a7}
\{ \vue \}_{\ep > 0} \ \mbox{bounded in} \ L^2(0,T ; W^{1,2}(\Omega; \R^3)).
\eF

In accordance with (\ref{a2}), (\ref{a3}), we may assume that
\bFormula{a8}
[r_\ep]_{\rm ess} \to r \ \mbox{weakly-(*) in}\ L^\infty(0,T; L^2(\Omega)),
\eF
and, taking (\ref{a4}) into account,
\bFormula{a8a}
[r_\ep]_{\rm res} \to 0 \ \mbox{in} \ L^\infty(0,T; L^q(\Omega))
\ \mbox{for any}\ 1 \leq q < \min\{ \gamma, 2 \}.
\eF
Moreover, by virtue of (\ref{a7}),
\bFormula{a9}
\vue \to \vc{U} \ \mbox{weakly in} \ L^2(0,T; W^{1,2}(\Omega; \R^3)),
\eF
passing to suitable subsequences as the case may be.

Letting $\ep \to 0$ in the weak formulation of the continuity equation (\ref{i11}), with
$b \equiv 0$, we obtain
\bFormula{a10}
\Div \vc{U} = 0 \ \mbox{a.a. in}\ (0,T) \times \Omega.
\eF

Finally, multiplying momentum balance (\ref{i12}) by $\ep$, we recover (\ref{i8})
\bFormula{a11}
\overline\vr\VC{-U_2}{ U_1} = p'(\overline\vr)\nabla_h r, \quad  \partial_{3} r = 0,
\eF
in particular, $r=r(x_1,x_2)$ is independent of the vertical variable, and
\bFormula{a12}
\vc{U}_h = \vc{U}_h(x_1,x_2), \
{\rm div}_h \vc{U}_h = 0,
\eF
which, together with (\ref{a10}), implies $U_3$ is independent of $x_3$. However, as~$\vc{U}$ satisfies the complete-slip boundary conditions (\ref{i5}) on $\partial \Omega$,
we may infer that
\bFormula{a13}
U_3 \equiv 0.
\eF

\section{Propagation of acoustic waves}
\label{aw}

Assume from now on, to simplify notation, that $p'(\Ov{\vr}) = 1$.
System~(\ref{i11}), (\ref{i12}) can be written in the form
\bFormula{aw1}
\ep \partial_t r_\ep + \Div \vc{V}_\ep = 0,
\eF
\bFormula{aw2}
\ep \partial_t \vc{V}_\ep + \left( \vc{g} \times \vc{V}_\ep + \Grad r_\ep\right) = \ep \vc{f}_\ep,
\eF
where we have set
\[
r_\ep = \frac{\vre - \Ov{\vr}}{\ep}, \ \vc{V}_\ep = \vre \vue,
\]
and
\[
\vc{f}_\ep = \Div \tn{S}(\Grad \vue) - \Div (\vre \vue \otimes \vue) - \frac{1}{\ep^2} \Grad \Big( p(\vre) - p'(\Ov{\vr})(\vre - \Ov{\vr}) - p(\Ov{\vr}) \Big).
\]
More precisely, system (\ref{aw1}), (\ref{aw2}) should be understood in the weak sense:
\bFormula{aw3}
\int_0^T \intO{ \Big( \ep r_\ep \partial_t \varphi + \vc{V}_\ep \cdot \Grad \varphi \Big) }
\ \dt = - \ep \intO{ r_{0,\ep} \varphi(0, \cdot) }
\eF
for any $\varphi \in \DC([0,T) \times \Ov{\Omega} )$,
\bFormula{aw4}
\int_0^T \intO{ \Big( \ep \vc{V}_\ep \cdot \partial_t \varphi - (\vc{g} \times \vc{V}_\ep) \cdot \varphi + r_\ep \Div \varphi \Big) } \ \dt = - \ep \int_0^T < \vc{f}_\ep, \varphi > \ \dt
\eF
\[
- \ep \intO{ \vr_{0,\ep} \vu_{0,\ep} \cdot \varphi(0, \cdot) },
\]
for any test function $\varphi \in \DC([0,T) \times \Ov{\Omega}; \R^3)$, $\varphi \cdot \vc{n}|_{\partial \Omega} = 0$, where
\[
- < \vc{f}_\ep, \varphi >
 = \int_{\Omega} \Big( \tn{S}(\Grad \vue) : \Grad \varphi
- (\vre \vue \otimes \vue): \Grad \varphi
\]
\[
- \frac{1}{\ep^2} \Big(
p(\vre) - p'(\Ov{\vr}) (\vre - \Ov{\vr}) - p(\Ov{\vr}) \Big) \Div \varphi \Big) \ \dx.
\]

It follows from the uniform bounds established in (\ref{a1} - \ref{a7}) that
\bFormula{aw5}
< \vc{f}_\ep , \varphi > = \intO{ \Big( \tn{F}^1_\ep : \Grad \varphi + \tn{F}^2_\ep :
\Grad \varphi \Big) },
\eF
with
\bFormula{aw6}
\{ \tn{F}^1_\ep \}_{\ep > 0} \ \mbox{bounded in}\ L^\infty(0,T; L^1(\Omega; \R^{3 \times 3})), \
\eF
\bFormula{aw7}
\{ \tn{F}^2_\ep \}_{\ep > 0} \ \mbox{bounded in}\ L^2(0,T; L^2(\Omega; \R^{3 \times 3})).
\eF

\subsection{Point spectrum of the acoustic propagator}

Consider an operator ${\mathcal B}$ defined, formally, in $L^2(\Omega) \times L^2(\Omega; \R^3)$,
\[
{\mathcal B} \VC{ r }{\vc{V}} \equiv \VC{ \Div \vc{V} }{ \vc{g} \times \vc{V} + \Grad r }.
\]
As a matter of fact, it is more convenient to work in the frequency space, meaning, we associate to a function $v$ its Fourier transform $\tilde v$
\[
\tilde v = \tilde v(\xi_h, k),\
\xi_h \equiv (\xi_1, \xi_2) \in \R^2, \ k \in Z,
\]
where
\[
\tilde v (\xi_h, k ) = \int_0^1 \int_{\R^2} \exp\Big(- {\rm i}(\xi_h \cdot x_h ) \Big) v(x_h, x_3) \ {\rm d}x_h \exp( - {\rm i}k x_3) \ {\rm d}x_3.
\]

%\subsubsection{Point spectrum of the operator ${\mathcal B}$}

We investigate the point spectrum of ${\mathcal B}$, meaning, we look for solutions of the eigenvalue problem
\bFormula{aw8}
\Div \vc{V} = \lambda r, \ \Grad r + \vc{g} \times \vc{V} = \lambda \vc{V},
\eF
or, in the Fourier variables,
$$
 {\rm i} \Big( \sum_{j=1}^2 \xi_j \tilde V_j + k \tilde V_3 \Big) - \lambda \tilde r = 0,
\  {\rm i} [\xi_1, \xi_2, k ] \tilde r - [  \tilde V_2, -\tilde V_1, 0 ] -
\lambda \tilde {\vc{V}} = 0.
$$

After a bit tedious but straightforward manipulation, we obtain
\bFormula{aw9}
\lambda^2 = - \mu, \ \mu = \frac{ 1 + |\xi|^2 + k^2 \pm \sqrt{ (1 + |\xi|^2 + k^2)^2 - 4 k^2} }{2};
\eF
whence the only eigenvalue is $\lambda = 0$, for which $k = 0$, and consequently, the space of eigenvectors coincides with the null-space of ${\mathcal B}$,
\bFormula{aw10}
{\rm Ker}({\mathcal B}) = \Big\{ [r, \vc{V}] \ \Big| \ r = r(x_1,x_2), \
\eF
\[
\vc{V} = [V_1(x_1,x_2),
V_2(x_1,x_2), V_3(x_1,x_2)] , \ {\rm div}_h  \vc{V}_h = 0, \ \nabla_h r = [V_2, -V_1] \Big\} .
\]

\subsection{RAGE theorem}

Our goal is to show that the component of the field $[r_\ep, \vc{V}_\ep]$, orthogonal to the null space ${\rm Ker}({\mathcal B})$ decays to zero on any compact subset of $\Omega$. To this end, we use the celebrated RAGE theorem in the following form (see Cycon et al.
\cite[Theorem 5.8]{CyFrKiSi}):
\bTheorem{aw1}
Let $H$ be a Hilbert space, $A: {\mathcal D}(A) \subset H \to H$ a self-adjoint operator, $C: H \to H$ a compact operator, and $P_c$ the orthogonal projection onto $H_c$, where
\[
H = H_c \oplus {\rm cl}_H \Big\{ {\rm span} \{ w \in H \ | \ w \ \mbox{an eigenvector of} \ A \} \Big\}.
\]

Then
\[
\left\| \frac{1}{\tau} \int_0^\tau \exp(-{\rm i} tA ) C P_c \exp(  {\rm i} tA ) \ \dt \right\|_{{\mathcal L}(H)} \to 0
\ \mbox{for}\ \tau \to \infty.
\]
\eT

In addition to the hypotheses of Theorem \ref{Taw1}, suppose that $C$ is non-negative and self-adjoint in $H$. Thus we may write
\[
\frac 1T \int_0^T \left< \exp\left( -{\rm i} \frac{t}{\ep} A \right) C
\exp \left( {\rm i} \frac{t}{\ep} A \right) P_c X, Y \right>_H {\rm d}t\leq  h(\ep) \| X \|_H \| Y \|_H,
\]
where $h(\ep) \to 0$ as $\ep \to 0$. Taking $Y = P_c X$ we deduce
\bFormula{aw11}
\frac 1T \int_0^T \left\| \sqrt{C} \exp \left( {\rm i} \frac{t}{\ep} A \right) P_c X \right\|_H^2 {\rm d}t
\leq  h(\ep) \| X \|^2_H.
\eF

Similarly, for $X \in L^2(0,T; H)$, we have
\begin{eqnarray}
\label {aw12}
& \displaystyle \frac 1{T^2}  \left\| \sqrt{C} P_c \int_0^t \exp \left( {\rm i} \frac{t-s}{\ep} A \right) X(s) \ {\rm d}s \right \|_{L^2(0,T; H)}^2  \\
 \leq & \displaystyle\frac 1T \int_0^T \int_0^T \left\| \sqrt{C} \exp \left( {\rm i} \frac{t-s}{\ep} A \right) P_c X(s)
\right\|^2_H \ {\rm d}t \ {\rm d}s
 \nonumber \\
 \leq & \displaystyle h(\ep) \int_0^T \left\| \exp \left( - {\rm i} \frac{s}{\ep} A \right) X(s)
\right\|^2_H \ {\rm d}s =  h(\ep) \int_0^T \| X(s) \|^2 \ {\rm d}s. \nonumber
\end{eqnarray}

%
%\bFormula{aw12}
%\left\| \sqrt{C} P_c \int_0^t \exp \left( {\rm i} \frac{t-s}{\ep} A \right) X(s) \ {\rm d}s \right \|_{L^2(0,T; H)}^2  \leq
%\eF
%\[
%T \int_0^T \int_0^T \left\| \sqrt{C} \exp \left( {\rm i} \frac{t-s}{\ep} A \right) P_c X(s)
%\right\|^2_H \ {\rm d}t \ {\rm d}s
%\]
%\[
%\leq T^2 h(\ep) \int_0^T \left\| \exp \left( - {\rm i} \frac{s}{\ep} A \right) X(s)
%\right\|^2_H \ {\rm d}s = T^2h(\ep) \int_0^T \| X(s) \|^2 \ {\rm d}s.
%\]

\subsection{Application of RAGE theorem}

For a fixed $M > 0$, we introduce a Hilbert space
\[
H = H_M \equiv \{ [r, \vc{V}] \ | \ \tilde r(\xi_h,k) = 0, \tilde
{\vc{V}} (\xi_h, k) = 0 \ \mbox{whenever}\ |\xi_h | + |k| > M \}.
\]
Let
\[
P_M : L^2(\Omega) \times L^2(\Omega; \R^3) \to H_M
\]
denote the associated orthogonal projection onto $H_M$.

Our goal is to apply RAGE theorem to the operators
\[
A = {\rm i}{\mathcal B}, \ C[v] = P_M[\chi v] ,\ \chi \in
\DC(\Omega), \ 0 \leq \chi \leq 1,
\]
considered on the Hilbert space $H_M$.

Going back to system (\ref{aw3}), (\ref{aw4}), we obtain that
\bFormula{aw13}
\ep \frac{{\rm d}}{{\rm d}t} \VC{r_{\ep,M} }{\vc{V}_{\ep,M}} +
{\mathcal B} \VC{r_{\ep, M} }{\vc{V}_{\ep,M}} = \ep
\VC{0}{\vc{f}_{\ep,M}},
\eF
where
\[
[r_{\ep,M}, \vc{V}_{\ep,M}] = P_M [r_\ep, \vc{V}_\ep],
\]
and
\[
\VC{0}{\vc{f}_{\ep,M}} \in H_M^* \approx H_M,
\]
\[
\left< \VC{0}{\vc{f}_{\ep,M}} , \VC{s}{\vc{w}} \right>_{H_M} = -
\intO{ \Big( \tn{F}^1_\ep : \Grad \vc{w} + \tn{F}^2_\ep : \Grad
\vc{w} \Big)}
\]
whenever $(s,w)\in H_M$.
Since
\[
\| \vc{w} \|_{W^{m,\infty} \cap W^{m,2}(\Omega; \R^3)} \leq c(m)
\| \vc{w} \|_{W^{m + 2, 2}(\Omega; \R^3)} \leq c M^{m+2} \| \vc{w} \|_{L^2(\Omega; \R^3)},
\]
we may use the uniform bounds (\ref{aw6}), (\ref{aw7}) in order to conclude that
\[
\left\| \VC{0}{\vc{f}_{\ep,M}} \right\|_{L^2(0,T; H_M)} \leq c(M)
\]
uniformly for $\ep \to 0$.

Writing solutions to (\ref{aw13}) by means of Duhamel's formula we get
\bFormula{aw14}
\VC{r_{\ep,M}}{\vc{V}_{\ep,M}} = \exp ({\rm i} A \frac{t}{\ep})
\VC{r_{\ep,M}(0)}{\vc{V}_{\ep,M}(0) } + \int_0^t \exp \left({\rm i}
\frac{t-s}{\ep} A \right) \VC{0}{\vc{f}_{\ep,M}} \ {\rm d}s;
\eF
whence a direct application of (\ref{aw11}), (\ref{aw12}), recalling that the only point spectrum is reduced to 0, yields
\bFormula{aw15}
Q^\perp \VC{ r_{\ep,M}}{\vc{V}_{\ep, M} } \to 0
\ \mbox{in} \ L^2((0,T) \times K; \R^4)) \ \mbox{as}\ \ep \to 0,
\eF
for any compact $K \subset \Ov{\Omega}$ and any fixed $M$, where we have denoted
\[
Q : L^2(\Omega) \times L^2(\Omega; \R^3)  \to {\rm Ker}({\mathcal B})
\]
the orthogonal projection onto the null space of ${\mathcal B}$.
Indeed observe that
\[
\left\| \sqrt{C} Q^\perp \VC{r_{\ep,M}}{\vc{V}_{\ep,M}}
\right\|^2_{H_M} = \left< C Q^\perp \VC{r_{\ep,M}}{\vc{V}_{\ep,M}},
Q^\perp \VC{r_{\ep,M}}{\vc{V}_{\ep,M}} \right>_{H_M}
\]
\[
= \intO{ \chi \left| Q^\perp \VC{r_{\ep,M}}{\vc{V}_{\ep,M}}
\right|^2 },
\]
where we have used the fact that $P_M$ and $Q$ commute.

Finally, a direct inspection of (\ref{aw14}) yields
\bFormula{aw16}
Q \VC{ r_{\ep,M}}{\vc{V}_{\ep, M} } \to \VC{ r_M }{\Ov{\vr} \vc{U}_M }
\ \mbox{in} \ L^2((0,T) \times K; \R^4)) \ \mbox{as}\ \ep \to 0,
\eF
where $r$ and $\vc{U}$ are the asymptotic limits identified through (\ref{a8} - \ref{a11}).

\subsection{Strong convergence of the velocity fields}

Relations (\ref{aw15}), (\ref{aw16}), together with (\ref{a7} - \ref{a9}), may be used
to obtain the desired conclusion
\bFormula{aw17}
\vue \to \vc{U} \ \mbox{in}\ L^2((0,T) \times K; \R^3)
\ \mbox{for any compact}\ K \subset \Omega.
\eF
Indeed, by virtue of (\ref{a8}), (\ref{a8a}), (\ref{aw15}), (\ref{aw16}), we obtain
\[
P_M [ \vue ] \to P_M [\vc{U}] \ \mbox{in} \ L^2((0,T) \times K; \R^3)
\]
for any fixed $M$, which, together with (\ref{a9}) and \emph{compactness} of the
embedding $W^{1,2}(K) \hookrightarrow L^2(K)$, yields (\ref{aw17}).

\section{The limit system}
\label{c}

\subsection{Identifying the limit system}

With the convergence established in (\ref{a8} - \ref{a9}), and (\ref{aw17}), it is not difficult to pass to the limit in the weak formulation (\ref{i11}), (\ref{i12}). To this end, we take
\[
\varphi \equiv  [\nabla^\perp_h \psi, 0], \ \psi \in \DC([0,T) \times \Omega)
\]
as a test function in momentum equation (\ref{i12}) to obtain
\bFormula{c1}
\int_0^T \intO{ \Big( \vre \vue \cdot \partial_t \varphi + \vre \vue \otimes \vue : \Grad \varphi - \frac{1}{\ep} \vre [\vue]_h \cdot \Grad \psi \Big) } \ \dt
\eF
\[
= - \intO{ \vr_{0,\ep} \vu_{0,\ep} \cdot \varphi (0, \cdot) } + \int_0^T \intO{
\tn{S}(\Grad \vue) : \Grad \varphi } \ \dt.
\]
Moreover, (\ref{aw3}) yields
\bFormula{c2}
\int_0^T \intO{ \Big( r_\ep \partial_t \psi + \frac{1}{\ep} \vre [\vue]_h \cdot \Grad \psi \Big) } = - \intO{ r_{0,\ep} \psi(0, \cdot) } .
\eF

Letting $\ep \to 0$ in (\ref{c1}), (\ref{c2}) we may infer that
 \[
\int_0^T \intO{ \Big( \Ov{\vr} \vc{U}_h \cdot \partial_t \nabla_h^\perp \psi +
\Ov{\vr} [ \vc{U}_h \otimes \vc{U}_h ] : \Grad ( \nabla^{\perp}_h \psi ) + r \partial_t \psi
\Big)}
\]
\[
= - \intO{ \Big( \Ov{\vr} \vc{U}_{0,h} \cdot \nabla^\perp_h \psi(0,\cdot)  +
r_0 \psi (0,\cdot) \Big)}
\]
\[
+ \int_0^T \intO{ \mu \nabla_h \vc{U}_h : \nabla (\nabla^\perp_h \psi ) } \ \dt.
\]

Moreover, as the limit functions are independent of $x_3$, we get,
\bFormula{c3}
\int_0^T \int_{\R^2} \Big( \Ov{\vr} \vc{U}_h \cdot \partial_t \nabla_h^\perp \psi +
\Ov{\vr} [ \vc{U}_h \otimes \vc{U}_h ] : \nabla_h ( \nabla^{\perp}_h \psi ) + r \partial_t \psi
\Big) \ {\rm d} x_h \ \dt
\eF
\[
= - \int_{\Rm^2} \left( \Ov{\vr} \left( \int_0^1\vc{U}_{0,h} \ {\rm d}x_3 \right)  \cdot \nabla^\perp_h \psi(0,\cdot)  +
\left( \int_0^1 r_0 \ {\rm d}x_3 \right) \psi (0,\cdot) \right) \ {\rm d}x_h
\]
\[
+ \int_0^T \int_{\R^2} \mu \nabla_h \vc{U}_h : \nabla_h (\nabla^\perp_h \psi ) \ {\rm d}x_h \ \dt
\]
for all $\psi \in \DC([0,T) \times \R^2)$.

Finally,
by virtue of (\ref{i8}), $\vc{U}_h = \nabla_h^\perp r$, and (\ref{c3}) coincides with  a weak formulation of~(\ref{i10}), (\ref{idat}). We have completed the proof of the convergence result, up to a subsequence, of Theorem \ref{Tw1}.

\subsection{Uniqueness for the limit system}
In this final section we shall prove that the limit system has a unique solution
provided the initial data are more regular. In order to do so we shall simply write an energy-type estimate on the difference of two solutions, called~$r_1$ and~$r_2$, associated with two initial data~$\tilde r_1$ and~$\tilde r_2$. This will provide a stability estimate, whose immediate consequence will be a uniqueness result. Notice that the diagnostic equation~(\ref{i8}) implies that~$\tilde r$ should be taken in~$W^{1,2} (\R^2)$.

The limit system writes
$$
\partial_t  ( \Delta_h r -  r  ) +
\nabla^\perp_h r \cdot \nabla_h ( \Delta_h r ) = \frac{\mu}{\Ov{\vr}} \Delta^2_h r
$$
recalling that for simplicity we have chosen~$p'(\Ov{\vr}) = 1$.
Multiplying (formally) this equation by~$ \Delta_h r $ and integrating over~$\R^2$ yields
$$
  \frac d{dt}\Big( \|\Delta_h r\|_{L^2}^2  +  \| \nabla_h r\|_{L^2}^2\Big) +  \frac{\mu}{\Ov{\vr}} \|\nabla_h\Delta_h  r\|_{L^2}^2 =0,
$$
whence the estimate
$$
\longformule{
 \|\Delta_h r(t)\|_{L^2}^2   +  \| \nabla_h r(t)\|_{L^2}^2 +   \frac{2\mu}{\Ov{\vr}} \int_0^t
 \| \nabla_h\Delta_h r(t')\|_{L^2}^2 \: dt' }{=  \|\Delta_h \tilde r \|_{L^2}^2   +  \| \nabla_h \tilde r \|_{L^2}^2.}
 $$
 Now suppose~$r_1$ and~$r_2$ are two solutions as described above, and define~$\delta := r_1-r_2$. Then of course~$\delta$ satisfies
 $$
 \partial_t  ( \Delta_h \delta -   \delta ) +
 \nabla^\perp_h \delta \cdot \nabla_h ( \Delta_h r_2)  + \nabla^\perp_h r_1 \cdot \nabla_h ( \Delta_h \delta) = \frac{\mu}{\Ov{\vr}} \Delta^2_h \delta
 $$
 with initial data~$\delta^0 = \tilde r_1- \tilde r_2$. Writing a similar energy estimate to the one above yields formally
 $$
 \longformule{
  \frac d{dt}\Big( \|\Delta_h \delta\|_{L^2}^2  +  \| \nabla_h  \delta\|_{L^2}^2\Big) +  \frac{2\mu}{\Ov{\vr}} \|\nabla_h \Delta_h \delta\|_{L^2}^2 }{= - \int_{\R^2}  \nabla^\perp_h \delta \cdot \nabla_h ( \Delta_h r_1)  \Delta_h\delta \: dx.}
 $$
 Then we simply write, by H\"older's inequality followed by Gagliardo-Nirenberg's inequality
 \begin{eqnarray*}
&\displaystyle \Bigl|  \int_{\R^2}  \nabla^\perp_h \delta \cdot \nabla_h ( \Delta_h r_1)  \Delta_h\delta \: dx \Bigr|
\leq  \| \nabla^\perp_h \delta \|_{L^4}\|\nabla_h  \Delta_h r_1   \|_{L^2} \| \Delta_h\delta\|_{L^4} \\
& \: \quad\quad \quad \displaystyle\leq C\| \nabla_h \delta \|_{L^2}^\frac12  \| \Delta_h\delta\|_{L^2}^\frac12\|\nabla_h  \Delta_h r_1   \|_{L^2}
 \| \Delta_h\delta\|_{L^2}^\frac12\| \nabla_h\Delta_h\delta\|_{L^2}^\frac12.
 \end{eqnarray*}
 This implies that
$$\longformule{
 \Bigl|  \int_{\R^2}  \nabla^\perp_h \delta \cdot \nabla_h ( \Delta_h r_1)  \Delta_h\delta \: dx \Bigr|
\leq     \frac{ \mu}{\Ov{\vr}} \|\nabla_h  \Delta_h \delta   \|_{L^2}^2 +   \| \nabla_h \delta \|_{L^2}^2 }{+ C \sqrt \frac{ \Ov{\vr}}{\mu}  \| \Delta_h\delta\|_{L^2}^2\| \nabla_h\Delta_hr_1\|_{L^2}^2.}
$$
Finally Gronwall's inequality allows to obtain
\begin{eqnarray*}
& \quad \quad \quad \quad \quad \quad \displaystyle\|\Delta_h \delta(t)\|_{L^2}^2 +  \| \nabla_h \delta (t)\|_{L^2}^2  +  \frac{\mu}{\Ov{\vr}} \int_0^t
 \| \nabla_h\Delta_h \delta (t')\|_{L^2}^2 \: dt' \\
  &    \leq \displaystyle \Big(\|\Delta_h \delta^0 \|_{L^2}^2   +  \| \nabla_h \delta^0\|_{L^2}^2 \Big) \exp \Big( C \sqrt \frac{ \Ov{\vr}}{\mu} \int_0^t \| \nabla_h\Delta_hr_1(t')\|_{L^2}^2\: dt'+Ct\Big).
 \end{eqnarray*}
This allows to conclude to stability, hence uniqueness for the limit system (leaving the usual regularization procedure to make the above arguments rigorous to the reader) provided
the initial datum enjoys the extra regularity stated in Theorem \ref{Tw1}.

\def\ocirc#1{\ifmmode\setbox0=\hbox{$#1$}\dimen0=\ht0 \advance\dimen0
  by1pt\rlap{\hbox to\wd0{\hss\raise\dimen0
  \hbox{\hskip.2em$\scriptscriptstyle\circ$}\hss}}#1\else {\accent"17 #1}\fi}


\begin{thebibliography}{1}

\bibitem{bmn} A. Babin, A. Mahalov, and B. Nicolaenko, Global
regularity of 3D rotating Navier--Stokes equations for resonant
domains,  {\it  Indiana University Mathematics Journal}, {\bf 48} (1999),
pages~1133--1176.

\bibitem{BDG} D. Bresch, B. Desjardins, B., D.  G\'erard-Varet,  Rotating fluids in a cylinder, {\it   Discrete Contin. Dyn. Syst.} {\bf  11}  (2004),  no. 1, 47--82.

 \bibitem{cdgg} J.-Y. Chemin, B. Desjardins, I. Gallagher and
E. Grenier, {\sl Basics of Mathematical Geophysics}, {\sl Oxford
University Press}, 2006, xii+250 pages.

\bibitem{CyFrKiSi}
H.L. Cycon, R.G. Froese, W.~Kirsch, and B.~Simon.
\newblock {\em Schr{\" o}dinger operators: with applications to quantum
  mechanics and global geometry}.
\newblock Texts and monographs in physics, Springer-Verlag, Berlin,Heidelberg,
  1987.

\bibitem{desjardinsgrenier} B. Desjardins and
E. Grenier,  Low Mach
number limit of compressible flows in the whole
space,
{\it Proceedings of the Royal Society of London}, {\bf 455},  pages
2271--2279, 1999.

\bibitem{dglm} B. Desjardins, E. Grenier, P.-L. Lions and N. Masmoudi, Incompressible limit for solutions of the isentropic Navier-Stokes equations with Dirichlet boundary conditions.  {\it J. Math. Pures Appl.} (9)  {\bf 78}  (1999),  no. 5, 461--471.

\bibitem{GSR2} I. Gallagher  \&  L. Saint-Raymond,  On the influence of the
   Earth's rotation on geophysical flows,  {\it Handbook of
   Mathematical  Fluid Dynamics}, {S. Friedlander and D. Serre Editors} Vol 4,  Chapter 5, 201-329, 2007.

\bibitem{FNP}
E.~Feireisl, A.~Novotn{\' y}, and H.~Petzeltov{\' a}.
\newblock On the existence of globally defined weak solutions to the
  {N}avier-{S}tokes equations of compressible isentropic fluids.
\newblock {\em J. Math. Fluid Mech.}, {\bf 3}:358--392, 2001.

\bibitem{klainermanmajda} S.
Klainerman and A.  Majda, Singular limits of quasilinear
hyperbolic systems with large parameters, and the
incompressible limit of compressible fluids, {\it
Communications on Pure and Applied Mathematics}, {\bf 34},
pages~481--524, 1981.


\bibitem{LI4}
P.-L. Lions.
\newblock {\em Mathematical topics in fluid dynamics, Vol.2, Compressible
  models}.
\newblock Oxford Science Publication, Oxford, 1998.

\bibitem{ukai} S. Ukai,  The incompressible limit
and the initial
layer of the compressible Euler
equation. {\it  Journal of  Mathematics of Kyoto
University} {\bf 26}, no. $2$,
pages~323--331, 1986.

\end{thebibliography}
 \end{document}